\documentclass[12pt]{article}
\usepackage{epsfig}
\usepackage{amssymb}

\def\P{{\mathbb P}}
\newtheorem{theorem}{Theorem}[section]
\newtheorem{lemma}[theorem]{Lemma}
\newtheorem{corollary}[theorem]{Corollary}
\newtheorem{proposition}[theorem]{Proposition}
\newtheorem{definition}[theorem]{Definition}

\title{Survival of inhomogeneous Galton-Watson processes}

\author{Erik Broman\footnote{Chalmers University of Technology, {\tt broman@math.chalmers.se}}
\and Ronald Meester\footnote{VU University Amsterdam,
Dept.\ of mathematics, De Boelelaan 1081a,
1081 HV Amsterdam, The Netherlands, {\tt rmeester@few.vu.nl}}}
\date{June 5, 2008}

\begin{document}

\maketitle

\begin{abstract}

We study survival properties of inhomogeneous Galton-Watson
processes. We determine the so-called branching number
(which is the reciprocal of the critical value for percolation) for these random
trees (conditioned on being infinite), which turns out to be an a.s.\ constant.
We also shed some light on the way the survival probability varies between the
generations. When we perform independent percolation on the family tree of an inhomogeneous
Galton-Watson process, the result is essentially a family of inhomogeneous Galton-Watson processes,
parametrized by the retention probability $p$. We provide growth rates, uniformly in $p$, of the percolation clusters, and also show uniform convergence of the survival probability from the $n$-th
level along subsequences. These results also establish, as a corollary, the supercritical continuity
of the percolation function. Some of our results are generalisations of
results by Lyons (1992).

\end{abstract}


\section{Introduction and main results}

We start by defining the main object of study in this paper, namely inhomogeneous Galton-Watson
processes.
Start with a root $o$ and let $L_1$ be the distribution of the (random)
number of offspring of the root. Proceed by letting each child (if any) of the root have
an i.i.d.\ number of offspring with distribution $L_2,$ and also let these offsprings be
independent of the number of children of the root. Given a sequence
$\{L_n\}_{n=1}^{\infty},$ we let $L_n$ be the offspring distribution of every individual
of generation $n-1.$ Sometimes we will treat $L_n$ as a random variable rather than as a
distribution,
this is standard abuse of notation.
The root is considered to be generation 0. Observe that if the distributions
$\{L_n\}_{n=1}^{\infty}$
all are the same we get a regular Galton-Watson process. Observe also that if
${\mathbb P}(L_n=l_n)=1$ for every $n$ and some sequence of numbers $\{l_n\}_{n=1}^{\infty},$
then we a.s.\ get a (deterministic) spherically symmetric tree, that is, a rooted tree in which any
two
vertices in the same generation have the same degree.

We denote the random family tree of such an inhomogeneous Galton-Watson process by $T$.
We will let $\overline{T}$ be a tree
with distribution equal to $T$ conditioned on survival, and we will also let $I \subset
\overline{T}$ be the
tree that consists of
those vertices $x \in \overline{T}$ (and the edges between them) that have infinitely many
descendents in
$\overline{T}$.
We will denote by $T_n$, $\overline{T}_n$ and $I_n$ the number of points in the $n$th generation
of $T$, $\overline{T}$ and $I$ respectively.

It is well known (see e.g.\ \cite{Ly2}) and not hard to see that $I$ is itself the family tree of an inhomogeneous Galton-Watson process; we will use this fact later on.

For inhomogeneous Galton-Watson processes we define
the survival probability $\theta_n$ from the $n$th generation, that is,
\[
\theta_n:=\lim_{m \rightarrow \infty}{\mathbb P}(T_m>0 \Big{|} T_n=1).
\]

For an infinite tree, a {\em cutset} $\pi$ is defined to be a finite set of edges such that every
infinite
path starting at the origin must contain at least one edge of the cutset. We denote
by $\Pi$ the set of all such cutsets. Any infinite tree $S$ has a so-called {\em branching number}
which is defined
as follows.

\begin{definition} \label{def1}
The {\em branching number} of an infinite tree $S$ with root $o$ is denoted by
${\rm br} S$ and defined by
\[
{\rm br} S:=\sup\left\{\lambda ; \inf_{\pi \in \Pi} \sum_{e \in \pi} \lambda^{-|e|}>0 \right\}.
\]
\end{definition}

The branching number is a very important property for trees (see \cite{R}).
For instance it is known that (see \cite{Lyons90}) the critical density $p_c(S)$ for independent
percolation (we are assuming that the reader is familiar with the concept of percolation,
otherwise please see \cite{G}) on $S$ is the reciprocal of the
branching number, that is,
$$
p_c(S)= 1/{\rm br}S.
$$
Closely related to the branching number is the {\em lower growth number} $\underline{\rm gr} S$
which
is
defined by
\[
\underline{\rm gr} S:=\liminf_{n \rightarrow \infty} S_n^{1/n},
\]
where $S_n$ denotes the number of vertices in the $n$th generation of $S$.
It is not hard to see that we always have ${\rm br} S \leq \underline{\rm gr} S,$
while equality is not always true. It is however well-known that if $S$ is spherically symmetric, then
${\rm br} S = \underline{\rm gr} S$.

We start with the following simple survival criterion. This result is essentially contained in
Proposition 4.15 of \cite{Ly2}, but we do give a different proof based on even earlier work in \cite{Agresti}. The reason is that some of the
elements in the proof will be used again later in this paper.

\begin{proposition} \label{lemma4}
For any inhomogeneous Galton-Watson process with offspring distributions $\{L_n\}_{n=1}^{\infty}$
we have that
\[
\liminf_{n \rightarrow \infty}({\mathbb E}[T_n])^{1/n}<1 \Rightarrow
\lim_{n \rightarrow \infty}{\mathbb P}(T_n>0)=0.
\]
Furthermore, if
\begin{equation} \label{eqn11}
\sup_n {\mathbb E}[L^2_n]=C_1< \infty
\end{equation}
and
\begin{equation} \label{eqn12}
\inf_n {\mathbb E}[L_n]=C_2>0,
\end{equation}
then we also have that
\[
\liminf_{n \rightarrow \infty}({\mathbb E}[T_n])^{1/n}>1 \Rightarrow
\lim_{n \rightarrow \infty}{\mathbb P}(T_n>0)>0.
\]
\end{proposition}

Next we have a result concerning the branching number of $\overline{T}$. A priori this is a random
variable, but it turns out that ${\rm br} \overline{T}$ is an almost sure
constant (under mild conditions).

\begin{theorem} \label{thm4}
Consider an inhomogeneous Galton-Watson process with offspring
distributions $\{L_n\}_{n=1}^{\infty}$ satisfying
(\ref{eqn11})
and (\ref{eqn12}). Assume also that
\begin{equation} \label{eqn13}
\liminf_{n \rightarrow \infty}{\mathbb E}[T_{n}]^{1/n}>1.
\end{equation}
Then we have that ${\rm br} \overline{T} = \liminf_{n \rightarrow \infty} {\mathbb E}[T_n]^{1/n}$,
$[\overline{T}]$-a.s.
\end{theorem}

\noindent
We make some remarks about this result.
\begin{enumerate}
\item In \cite{Ly2}, it is proved that a.s.,
\[
{\rm br}\overline{T}=\liminf_{n \rightarrow \infty} {\mathbb E}[T_n]^{1/n}
\]
under the assumption that $\sup_{n}||L_n||_{\infty}<\infty$. It is claimed in \cite{Ly2} that this
assumption cannot be weakened much; our results show that if one adds the very natural condition
(\ref{eqn13}), then in fact one can significantly weaken the assumptions.

\item
Naively one might believe that this result would follow from easy arguments.
For instance one might try the following approach:
Define a new inhomogeneous Galton-Watson tree $T'$ by performing
percolation on $T$
with probability for being open equal to $p$. Depending on whether $p$ was smaller or
greater than $1/(\liminf E[T_n]^{1/n})$, we get from Proposition \ref{lemma4} that
$T'$ dies out a.s./survives with positive probability (respectively), concluding the
argument. However one then misses the point that
the fact that $T'$ survives with positive probability if $p>1/(\liminf
E[T_n]^{1/n})$ does not lead to the conclusion that ${\rm brT} \geq \liminf
E[T_n]^{1/n}$. Indeed, it is imaginable that with positive probability
${\rm brT}=\liminf E[T_n]^{1/n}-\delta$ for some positive $\delta$ and with positive
probability ${\rm brT}=\liminf E[T_n]^{1/n}.$ If this were true, $T$' would still
survive with positive probability for the indicated $p$.

\end{enumerate}

The following result about the behaviour of $\theta_n$ will be needed in the proof of Theorem
\ref{thm4} but is also quite interesting in its own right. It is not to be expected that $\theta_n$ is in
general bounded
away from 0 since one can always insert any finite number of generations of degree one in the tree.
 However, it is the case
that
there is a subsequence along which $\theta_n$ is bounded away from 0.

\begin{proposition} \label{lemma5}
Consider an inhomogeneous Galton-Watson process with offspring distributions
$\{L_n\}_{n=1}^{\infty}$
satisfying (\ref{eqn11}), (\ref{eqn12}) and (\ref{eqn13}).
Then there exists a sequence $\{n_k\}_{k=1}^{\infty}$
of increasing integers and a constant $C>0$ such that for all $k \geq 1$,
\[
\theta_{n_k} \geq C.
\]
\end{proposition}

Next, we study bond percolation on $I$. Note that
$p_c(I)=p_c(\overline{T})$, since pruning a tree does not change its critical probability.
We already noted that $I$ itself is the family tree of an inhomogeneous Galton-Watson process, and
when we perform independent bond percolation on $I$, the resulting component of the origin, to be
denoted by $I^p$, also constitutes a family tree of an inhomogeneous Galton-Watson process. Therefore, general results about inhomogeneous Galton-Watson
processes automatically apply to percolation on $I$. However, being equipped with a parameter $p$ now,
we will derive survival estimates {\em uniformly} in $p$.
We remark that a special case of inhomogeneous Galton-Watson processes results
from starting with a deterministic spherically symmetric tree and performing percolation on that
tree. One more piece of notation: the number of vertices in $I^p$
at distance $n$ from the root is denoted by $I_n^p$. Also, in this paper, we use various coupling constructions. To facilitate this,
all processes,
for all values of $p$, are jointly constructed in the obvious way. Consequently, as in the
previous example, we
will express the
$p$-dependence in the events rather than in the measure.

In light of Theorem \ref{thm4}, one might expect that for any
$\epsilon>0,$
\[
\lim_{n \rightarrow \infty}{\mathbb P}(0<I_n^p<((1-\epsilon){\rm br} \overline{I^p} )^n)=0.
\]
\noindent
In fact, we have the next, much stronger statement.

\begin{theorem} \label{propunifconv}
Consider an inhomogeneous Galton-Watson process satisfying (\ref{eqn11}),
(\ref{eqn12}) and (\ref{eqn13}), with family tree $T$, and let $\epsilon>0.$
For $p_c(\overline{T}) <p_1 \leq 1$ it is the case that
\[
\lim_{n \rightarrow \infty}{\mathbb P}(0<I_n^p<((1-\epsilon){\rm br}\overline{I^{p_1}})^n)=0,
\]
uniformly in $p \in [p_1,1].$
\end{theorem}

Note that the pointwise (in $p$) convergence in Theorem \ref{propunifconv} is almost a triviality. The
whole point of the theorem is proving the uniform convergence.

Theorem \ref{propunifconv} combined with Proposition \ref{lemma5} will in turn lead us to our
next result. Here we define
$$
\theta(p):={\mathbb P}(|I^p|=\infty).
$$

\begin{proposition} \label{thm2}
Consider an inhomogeneous Galton-Watson process satisfying (\ref{eqn11}), (\ref{eqn12}) and
(\ref{eqn13}),
and let $p_1>p_c(\overline{T}).$ Then there exists a
sequence of increasing integers $\{n_k\}_{k=1}^{\infty}$ such that
\[
\theta(p)=\lim_{k \to \infty} \P(I_{n_k}^p >0),
\]
uniformly on $[p_1,1].$
\end{proposition}

This result also leads to continuity of the percolation function above $p_c$
for random trees.

\begin{corollary} \label{corr1}
Consider an inhomogeneous Galton-Watson process satisfying (\ref{eqn11}), (\ref{eqn12}) and
(\ref{eqn13}).
Then the percolation function $\theta(p)$ is continuous above $p_c(\overline{T}).$
In particular, on any spherically symmetric tree $S$ with uniformly bounded degrees, the percolation
function is continuous above $p_c(S)$.
\end{corollary}
In fact, one can also use Theorem \ref{thm4} to construct a more or less classical proof of this
result.
As an interesting side remark, we mention that the route via Proposition \ref{thm2}
also has a counterpart on $\mathbb{Z}^d$, and gives a new proof for the continuity of the
percolation function
in that context. This proof does in fact give a rate of convergence for the natural approximations
of the
percolation function; we discuss these continuity matters in the last section.

In contrast to our last corollary, we have the following example of a tree
for which the percolation function is not continuous
above $p_c$. To
construct such a tree, we use a result in \cite{Ly2}, a special case of which says that there is
percolation with positive probability on a spherically symmetric tree $S$ with parameter $p$,
if and only if
$$
\sum_{n=1}^{\infty} \frac{p^{-n}}{S_n} < \infty.
$$
To construct an example, we first take a spherically symmetric tree $S$ which is such that $S_n$
is of the order $2^n n^2$. It follows from the above that $p_c(S)=1/2$ and that $\theta_{S}(1/2)>0$.
Next, we take a regular tree $S'$ with common degree 4. It is well-known that $p_c(S')=1/3$. We then
construct a tree $S''$ by joining the roots of $S$ and $S'$ by a single edge. It is easy to see
that $p_c(S'')=1/3$ and that $\theta_{S''}$ is discontinuous at $1/2$.

\medskip
Theorem \ref{thm4} along with Propositions \ref{lemma4} and \ref{lemma5} will be
proved in the next section.  All the other results are proved in Section \ref{sec3}. The
issues about continuity of the percolation function are discussed in Section \ref{sec4}.

\section{Proof of Proposition \ref{lemma4}, Theorem \ref{thm4} and Proposition \ref{lemma5}}
\label{secbr}
We start by defining a useful probability generating function by
\[
h(n,s):=\sum_{j=0}^{\infty}{\mathbb P}(L_n=j)s^j, \ \forall n \geq 1.
\]
It is known (see \cite{Agresti}) that if $h''(n,1)<\infty,$ for every $n,$
then
for all $n \geq 1$ we have
\begin{equation} \label{eqn10}
\left[{\mathbb E}[T_n]^{-1}+\sum_{j=1}^{n}\frac{h''(j,1)}{h'(j,1)} {\mathbb
E}[T_{j}]^{-1}\right]^{-1}
\leq {\mathbb P}(T_n>0).
\end{equation}
Of course we have
\[
h'(n,1)=\sum_{j=0}^{\infty}j{\mathbb P}(L_n=j)={\mathbb E}[L_n],
\]
and
\[
h''(n,1)=\sum_{j=0}^{\infty}j(j-1){\mathbb P}(L_n=j)={\mathbb E}[L^2_n]-{\mathbb E}[L_n].
\]

\noindent
We can now proceed with the proof of Proposition \ref{lemma4}.

\medskip\noindent
{\bf Proof of Proposition \ref{lemma4}.}
The proof of the first statement is easy.
Assume that $$\liminf_{n \rightarrow \infty}({\mathbb E}[T_n])^{1/n}=a<1,$$ then we get that
for any $\epsilon>0$ such that $a(1+\epsilon)<1,$ there exists a sequence $\{n_k\}_{k=1}^{\infty}$
such
that
\[
{\mathbb P}(T_{n_k}>0) \leq {\mathbb E}[T_{n_k}] \leq (a(1+\epsilon))^{n_k},
\]
so that
\[
\lim_{n \rightarrow \infty}{\mathbb P}(T_n>0)=0.
\]

For the second statement we start by observing that condition (\ref{eqn11}) gives us that
$h''(n,1)={\mathbb E}[L^2_n]-{\mathbb E}[L_n]<\infty$ for every $n.$ Of course this does not require
the
full statement of equation (\ref{eqn11}) which will be needed later.
In turn, this gives us that inequality (\ref{eqn10}) is valid for every $n,$ and therefore we need to
show
that
\begin{eqnarray} \label{eqn15}
\lefteqn{\limsup_{n \rightarrow \infty} \left[
{\mathbb E}[T_n]^{-1}+\sum_{j=1}^{n}\frac{h''(j,1)}{h'(j,1)} {\mathbb E}[T_{j}]^{-1}\right]^{-1}}
\\
& & =\limsup_{n \rightarrow \infty} \left[{\mathbb E}[T_n]^{-1}+\sum_{j=1}^{n}
\frac{{\mathbb E}[L^2_j]-{\mathbb E}[L_j]}{{\mathbb E}[L_j]} {\mathbb E}[T_{j}]^{-1}\right]^{-1}
>0. \nonumber
\end{eqnarray}
To this end, we observe that by equations (\ref{eqn11}) and (\ref{eqn12})
\[
\sup_{j}\frac{{\mathbb E}[L^2_j]-{\mathbb E}[L_j]}{{\mathbb E}[L_j]} \leq \frac{C_1}{C_2}=C<\infty.
\]

Since $\liminf_{n \rightarrow \infty}({\mathbb E}[T_n])^{1/n}>1$ there exists a constant $b>1$
and an $N$ such that
for all $n \geq N$,
\[
{\mathbb E}[T_n]>b^n.
\]
Therefore, for some constant $D<\infty,$
\begin{eqnarray*}
\lefteqn{
{\mathbb E}[T_n]^{-1}+\sum_{j=1}^{n}\frac{h''(j,1)}{h'(j,1)} {\mathbb E}[T_{j}]^{-1}} \\
& & \leq {\mathbb E}[T_n]^{-1}+C\sum_{j=1}^{n} {\mathbb E}[T_{j}]^{-1}
\leq D+C\sum_{j=N}^{\infty} b^{-j}< \infty.
\end{eqnarray*}
Since the right hand side of the above inequality is independent of $n,$
inequality (\ref{eqn15}) is valid and that concludes the proof.
\hfill{$\Box$}

\medskip\noindent
We continue by proving Proposition \ref{lemma5}.

\medskip \noindent{\bf Proof of Proposition \ref{lemma5}.}
Let $\{X_i\}_{i\geq 1}$ be i.i.d. with distribution according to $T_n$ conditioned
on the event that $T_{\ell}=1.$ Observe that for $n \geq \ell$
\[
T_n=\sum_{k=1}^{T_{\ell}}X_k,
\]
so that (using Wald's lemma)
\begin{equation} \label{eqn16}
{\mathbb E}[X_1]={\mathbb E}[T_n \Big{|} T_{\ell}=1]=\frac{{\mathbb E}[T_n]}{{\mathbb E}[T_{\ell}]}.
\end{equation}
Observe that by condition (\ref{eqn11}) we can use inequality (\ref{eqn10}) to conclude that
for $n \geq \ell,$
\begin{equation} \label{eqn19}
\left[{\mathbb E}[T_n\Big{|} T_{\ell}=1]^{-1}
+\sum_{j=\ell+1}^{n}\frac{h''(j,1)}{h'(j,1)} {\mathbb E}[T_{j}\Big{|} T_{\ell}=1]^{-1}\right]^{-1}
\leq {\mathbb P}(T_n>0\Big{|} T_{\ell}=1).
\end{equation}
We will show that there exists a sequence $\{n_k\}_{k=1}^{\infty}$
of increasing integers and a constant $C<\infty$ such that for all $k \geq 1$ and for all $n \geq
n_k$,
\begin{equation} \label{eqn13a}
{\mathbb E}[T_n\Big{|} T_{n_k}=1]^{-1}
+\sum_{j=n_k+1}^{n}\frac{h''(j,1)}{h'(j,1)} {\mathbb E}[T_{j}\Big{|} T_{n_k}=1]^{-1}\leq C.
\end{equation}
This will give us that
for all $k \geq 1$ we have
\[
\lim_{n \rightarrow \infty} {\mathbb P}(T_n>0\Big{|} T_{n_k}=1) \geq \frac{1}{C},
\]
proving the lemma. To that end, observe that as in the proof of
Lemma \ref{lemma4} there exists a constant $C_3$ such that
for $ n \geq \ell$,
\begin{eqnarray} \label{eqn18}
\lefteqn{{\mathbb E}[T_n\Big{|} T_{\ell}=1]^{-1}
+\sum_{j=\ell+1}^{n}\frac{h''(j,1)}{h'(j,1)} {\mathbb E}[T_{j}\Big{|} T_{\ell}=1]^{-1}} \\
& & \leq {\mathbb E}[T_n\Big{|} T_{\ell}=1]^{-1} \nonumber
+C_3\sum_{j=\ell+1}^{n} {\mathbb E}[T_{j}\Big{|} T_{\ell}=1]^{-1} \\
& & \leq (C_3+1)\sum_{j=\ell+1}^{n} {\mathbb E}[T_{j}\Big{|} T_{\ell}=1]^{-1} \nonumber
=(C_3+1){\mathbb E}[T_{\ell}]\sum_{j=\ell+1}^{n} \frac{1}{{\mathbb E}[T_{j}]},
\end{eqnarray}
where we use equation (\ref{eqn16}) in the last equality.
Therefore, showing that there exists a sequence $\{n_k\}_{k=1}^{\infty}$
of increasing integers and a constant $C<\infty$ such that for all $k$ we have
\[
{\mathbb E}[T_{n_k}]\sum_{j=n_k+1}^{\infty} \frac{1}{{\mathbb E}[T_{j}]} \leq C
\]
will give us equation (\ref{eqn13a}).

We divide the proof into three cases. First however, define
\[
m:=\liminf_{n \rightarrow \infty}{\mathbb E}[T_{n}]^{1/n}>1.
\]
In the first case, we have that
${\mathbb E}[T_{n}]^{1/n}<m$
for infinitely many $n.$  We can then conclude that there exists $n_1$, defined
to be the largest integer such that ${\mathbb E}[T_{n_1}]^{1/n_1}=\min_{n \geq 1} {\mathbb
E}[T_{n}]^{1/n}.$
Having defined $n_k,$ we can then define $n_{k+1}$
to be the largest integer greater than $n_k$ such that
${\mathbb E}[T_{n_{k+1}}]^{1/n_{k+1}}=\min_{n >n_k} {\mathbb E}[T_{n}]^{1/n}.$
Let $\epsilon_k$ be defined through ${\mathbb E}[T_{n_k}]^{1/n_k}=m(1-\epsilon_k).$
Observe that by definition of $n_k,$ ${\mathbb E}[T_{n}]^{1/n}\geq m (1-\epsilon_k)$
for every $n \geq n_k$ and also that $\epsilon_k >0$ for every  $k,$ and finally that
$\epsilon_k \rightarrow 0,$ as $ k \rightarrow \infty.$ Therefore,
\begin{eqnarray*}
\lefteqn{{\mathbb E}[T_{n_k}]\sum_{j=n_k+1}^{\infty} \frac{1}{{\mathbb E}[T_{j}]}
\leq {\mathbb E}[T_{n_k}]\sum_{j=n_k+1}^{\infty}\frac{1}{(m(1-\epsilon_k) )^{j}}} \\
& & = (m(1-\epsilon_k) )^{n_k}\sum_{j=1}^{\infty}\frac{1}{(m(1-\epsilon_k) )^{n_k+j}}
=\sum_{j=1}^{\infty}\frac{1}{(m(1-\epsilon_k) )^{j}}.
\end{eqnarray*}
There exists a $K$ such that $m(1-\epsilon_k)>1$ for every $k \geq K.$ For $k \geq K,$ the right
hand
side
of the above equation is then bounded by some constant $D_k< \infty.$ Furthermore, we can take
$D_k \geq D_{k+1}$ and conclude that for all $k \geq K$,
\[
{\mathbb E}[T_{n_k}]\sum_{j=n_k+1}^{\infty} \frac{1}{{\mathbb E}[T_{j}]}\leq D_K < \infty.
\]

For the second and third case, we have that ${\mathbb E}[T_{n}]^{1/n}<m$
for only finitely many $n.$ We can therefore find $N$ large enough so that
${\mathbb E}[T_{n}]^{1/n}\geq m$ for every $n \geq N.$ We have that for every $n,$
${\mathbb E}[T_{n}]^{1/n}=m(1+a(n)),$ where the sequence of numbers $\{a(n)\}_{n=1}^{\infty}$
is such that $a(n) \geq 0$ for every $n \geq N.$

The second case is if
$\liminf_{n \rightarrow \infty}(1+a(n))^n=C_4$ for some constant $C_4< \infty.$ Then there exists a
sequence of strictly increasing integers $\{n_k\}_{k=1}^{\infty}$ such that
$(1+a(n_k))^{n_k} \leq 2C_4$ for every $k \geq 1.$ By also requiring that $n_1 \geq N,$ we get that
\begin{eqnarray*}
\lefteqn{{\mathbb E}[T_{n_k}]\sum_{j=n_k+1}^{\infty} \frac{1}{{\mathbb E}[T_{j}]}
 \leq {\mathbb E}[T_{n_k}]\sum_{j=n_k+1}^{\infty} \frac{1}{m^{j}}} \\
& & =m^{n_k}(1+a(n_k))^{n_k}\sum_{j=1}^{\infty} \frac{1}{m^{n_k+j}}
 \leq 2C_4\sum_{j=1}^{\infty} \frac{1}{m^{j}} < \infty.
\end{eqnarray*}

The third case is if $\lim_{n \rightarrow \infty}(1+a(n))^n=\infty.$
We can then find a sequence $\{n_k\}_{k=1}^{\infty}$
(much as in the first case) such that for every $k \geq 1,$ $(1+a(n))^n \geq (1+a(n_k))^{n_k}$
for every $n \geq n_k.$ By again requiring that $n_1 \geq N,$ we get that
\begin{eqnarray*}
\lefteqn{{\mathbb E}[T_{n_k}]\sum_{j=n_k+1}^{\infty} \frac{1}{{\mathbb E}[T_{j}]}
\leq {\mathbb E}[T_{n_k}]\sum_{j=n_k+1}^{\infty}\frac{1}{(m(1+a(n_k)))^{j}}} \\
& & =\frac{(m(1+a(n_k)))^{n_k}}
{(m(1+a(n_k)))^{n_k}}\sum_{j=1}^{\infty}\frac{1}{(m
(1+a(n_k)))^{j}}
\leq \sum_{j=1}^{\infty}\frac{1}{m^{j}} < \infty.
\end{eqnarray*}

We can therefore conclude that
there exists a constant $C=C(\{L_n\}_{n=1}^{\infty})<\infty$
and a sequence of strictly increasing integers
$\{n_k\}_{k=1}^{\infty}$ such that
for all $k \geq 1$,
\[
{\mathbb E}[T_{n_k}]\sum_{j=n_k+1}^{\infty} \frac{1}{{\mathbb E}[T_{j}]}\leq C.
\]
This concludes the proof.
\hfill{$\Box$}

\medskip\noindent
We are now ready to prove Theorem \ref{thm4}.

\medskip\noindent
{\bf Proof of Theorem \ref{thm4}.}
Using that $p_c(I)^{-1}=p_c(\overline{T})^{-1}={\rm br \overline{T}}$, we need to show that
$$
p_c(I)^{-1}=
\liminf_{n \rightarrow \infty} ({\mathbb E}[T_n])^{1/n}.
$$
We will do this by
first proving that $p_c(I)^{-1}=
\liminf_{n \rightarrow \infty} (\theta_{n}{\mathbb E}[T_n])^{1/n}$ and then proving that
$$\
\liminf_{n \rightarrow \infty} (\theta_{n}{\mathbb E}[T_n])^{1/n}
=\liminf_{n \rightarrow \infty} {\mathbb E}[T_n]^{1/n}.
$$

Consider the offspring distribution $L'_1$ of the root of $I.$ Let $T^i$ be the tree
consisting of
child number $i\in\{1,\ldots,L_1\}$ of the root of $T$ and all the descendents of this child.
Define also $N_{1,\infty}=|\{T^i: |T^i|=\infty, i=1,\ldots,L_1\}|.$ It is not hard
to see that for $k \geq 1,$
\[
{\mathbb P}(L'_1=k)={\mathbb P}(N_{1,\infty}=k \Big{|} N_{1,\infty} \geq 1)
=\frac{{\mathbb P}(N_{1,\infty}=k)}{\theta}.
\]
Furthermore, letting $Y_i$ be i.i.d. ${\rm Bin}(1,\theta_1)$ random variables and using Wald's
lemma we get that
\[
{\mathbb E}[L'_1]=\frac{1}{\theta}{\mathbb E}[N_{1,\infty}]
=\frac{1}{\theta}{\mathbb E} \left[\sum_{i=1}^{L_1} Y_i\right]
=\frac{1}{\theta}{\mathbb E} \left[Y_1\right]{\mathbb E} \left[L_1\right]
=\frac{\theta_1}{\theta}{\mathbb E}[L_1].
\]
Furthermore, this argument holds for any generation $n$ and therefore we have for all $ n \geq 1$,
\begin{equation}
\label{generations}
{\mathbb E}[L'_n]=\frac{\theta_n}{\theta_{n-1}}{\mathbb E}[L_n].
\end{equation}

Now, perform independent percolation on $I$ with parameter $p,$ thus creating
a random graph that we denote by ${\cal I}^p.$ Recall that $I^p$ is the component of the
root of this graph.
Obviously, $I^p$ is the family tree of an inhomogeneous
Galton-Watson process with some offspring distributions
$\{L''_n\}_{n=1}^{\infty}.$ Furthermore, trivially
\[
{\mathbb E}[L''_n]=p{\mathbb E}[L'_n]=p\frac{\theta_n}{\theta_{n-1}}{\mathbb E}[L_n] \ \ \forall
n\geq
1.
\]

Recall that $I_n^p$ is the number of vertices in $I^p$ at distance $n$ from the root
 and recall that we defined $I_n$ similarly. We have, using a
standard result from the theory of branching processes and (\ref{generations}), that
\begin{equation}
{\mathbb E}[I_n^p]=p^n{\mathbb E}[I_n]
 =p^n \prod_{i=1}^{n}{\mathbb E}[L'_i]
=p^n \prod_{i=1}^{n}\frac{\theta_i}{\theta_{i-1}}{\mathbb E}[L_i]
=p^n \frac{\theta_{n}}{\theta}{\mathbb E}[T_n].
\end{equation}
Therefore,
\begin{equation} \label{eqn14}
\liminf_{n \rightarrow \infty}({\mathbb E}[I_n^p])^{1/n}
=p\liminf_{n \rightarrow \infty}(\frac{\theta_{n}}{\theta}{\mathbb E}[T_n])^{1/n}
=p\liminf_{n \rightarrow \infty}(\theta_{n}{\mathbb E}[T_n])^{1/n}.
\end{equation}
We would like to use Proposition \ref{lemma4} and Proposition \ref{lemma5} on $I^p.$ However,
before we can do that we need to show that the offspring distributions $\{L''_n\}_{n=1}^{\infty}$
satisfies conditions (\ref{eqn11}) and (\ref{eqn12}). When we use Proposition \ref{lemma5} we will
assume that condition (\ref{eqn13}) is satisfied; the details will become clear.

For some vertex $x$ in generation $n-1,$ let $T^i_x$ be the tree consisting of
child number $i\in\{1,\ldots,L_n\}$ of $x$ and all the descendents of this child.
Define $N_{n,\infty}=|\{T^i_x: |T^i_x|=\infty, i=1,\ldots,L_n\}|,$ and observe that
the distribution of this random variable is trivially independent of the specific choice of $x$ in
generation $n-1.$
Let $Y_i^n$ be i.i.d. ${\rm Bin}(1,\theta_n)$ and observe that
\begin{eqnarray*}
{\mathbb E}[(L''_n)^2] & \leq & {\mathbb E}[(L'_n)^2]=\sum_{j=1}^{\infty}j^2{\mathbb P}(L'_n=j) \\
&  = & \sum_{j=1}^{\infty}j^2{\mathbb P}(N_{n,\infty}=j \Big{|}N_{n,\infty} \geq 1 ) \\
& = & \frac{{\mathbb E}[N_{n,\infty}^2]}{\theta_{n-1}}
=\frac{{\mathbb E}\left[{\mathbb E}\left[\left(\sum_{i=1}^{L_n}Y_i^n\right)^2
\Big{|} L_n\right]\right]}{\theta_{n-1}} \\
& \leq & \frac{{\mathbb E}\left[{\mathbb E}\left[L_n \sum_{i=1}^{L_n}(Y_i^n)^2
\Big{|} L_n\right]\right]}{\theta_{n-1}} \\
& = & \frac{{\mathbb E}\left[L_n \sum_{i=1}^{L_n}{\mathbb E}\left[Y_i^n
\Big{|} L_n\right]\right]}{\theta_{n-1}}
=\frac{{\mathbb E}\left[L_n^2 \theta_n \right]}{\theta_{n-1}}.
\end{eqnarray*}
In the second inequality we use that for any real numbers $a_1, \ldots, a_n$ we have that
$(a_1+\cdots + a_n)^2 \leq n(a_1^2+\cdots+a_n^2).$
Obviously we must also have that
\[
\theta_{n-1} \geq {{\mathbb P}(L_n>0)\theta_{n}},
\]
and we can use Cauchy-Schwarz to see that
\[
{\mathbb E}[L_n]^2= {\mathbb E}[L_n I_{\{L_n>0\}}]^2\leq {\mathbb P}(L_n>0)  {\mathbb E}[L_n^2].
\]
Therefore,
\[
\frac{{\mathbb E}\left[L_n^2 \theta_n \right]}{\theta_{n-1}}
 \leq \frac{{\mathbb E}\left[L_n^2 \right]}{{\mathbb P}(L_n>0)}
\leq  {\mathbb E}\left[L_n^2 \right] \frac{{\mathbb E}\left[L_n^2 \right]}{{\mathbb E}\left[L_n
\right]^2}
\leq \frac{C_1^2}{C_2^2}<\infty.
\]
Furthermore
\[
\inf_{n}{\mathbb E}[L''_n]=p\inf_{n}{\mathbb E}[L'_n] \geq p,
\]
since ${\mathbb E}[L'_n] \geq 1$ for every $n.$

We can now proceed to use
Proposition \ref{lemma4} with equation (\ref{eqn14}) to see that $I^p$ survives with
positive probability if $p>(\liminf_{n \rightarrow \infty}(\theta_{n}{\mathbb E}[T_n])^{1/n})^{-1}$
while it dies out a.s. if  $p<(\liminf_{n \rightarrow \infty}(\theta_{n}{\mathbb
E}[T_n])^{1/n})^{-1}.$

This is not quite enough for our purposes: it could be the case that with positive probability, $I$
is
such that $I^p$ a.s. dies out. Since we want to make a statememt about almost all trees $I$,
we
argue
that in fact, if  $p>1/\liminf_{n \rightarrow \infty}(\theta_{n}{\mathbb E}[T_n])^{1/n}$, then
${\cal
I}^p$
contains an infinite component with probability 1 as our next argument shows.

Assume therefore that
$\liminf_{n \rightarrow \infty}({\mathbb E}[I^p_n])^{1/n}=
p\liminf_{n \rightarrow \infty}(\theta_n{\mathbb E}[T_n])^{1/n}>1.$ This is condition (\ref{eqn13})
for $I^p.$ Construct the tree $I^p$ by letting $I^p_1$ have distribution equal to $L_1''.$ Proceed
by letting $I^p_2$ be the sum $\sum_{i=1}^{I^p_1}L''_{2,i},$ where $\{L''_{2,i}\}_{i=1}^{\infty}$
are i.i.d.\ with distribution equal to $L_2''$ and let them also be independent of everything
else. Continuing in this fashion, we have two possibilities.
First we may find that
$I_n^p>0$ for every $n.$ Second, we might instead find that for some $n$, we have $I_n^p=0.$ If this
is
the
case, there exists some integer $n_{k_1} >n$ in the subsequence dictated by Proposition
\ref{lemma5}.
However, since $I$ is infinite we must have that ${\cal I}^p$ contains a subtree
(possibly consisting
of only one vertex )
with the root being some vertex at level $n_{k_1}.$ Construct this subtree in the same way as we
constructed $I^p$ above. This subtree has some probability to survive which is by Proposition
\ref{lemma5}
uniformly bounded away from 0. It is also easy to see that the event of survival
of this subtree is conditionally independent of the part of ${\cal I}^p$ examined so far (up to
generation $n$).

If again we find that this subtree is finite,
we continue in the same way. Since all
the subtrees that we pick have uniformly positive probability to survive by Proposition
\ref{lemma5} and the survival of them are conditionally independent we see that
${\cal I}^p$ must contain an infinite component with
probability 1.
We therefore conclude that
\[
{\mathbb P}({\cal I}^p \textrm{ has an infinite component})
=\left\{
\begin{array}{cc}
1, & p>1/\liminf_{n \rightarrow \infty}(\theta_{n}{\mathbb E}[T_n])^{1/n}, \\
0, & p < 1/\liminf_{n \rightarrow \infty}(\theta_{n}{\mathbb E}[T_n])^{1/n}.
\end{array}
\right.
\]
This is the same as saying that for almost every $I,$ we will after performing percolation with
parameter
$p$ on $I$
a.s.\ get an infinite component if $p>1/\liminf_{n \rightarrow \infty}(\theta_{n}{\mathbb
E}[T_n])^{1/n}$
while we will a.s. not get an infinite component if
$p<1/\liminf_{n \rightarrow \infty}(\theta_{n}{\mathbb E}[T_n])^{1/n}.$
It follows that for almost every $I$ the probability that the component of the root is infinite is
positive if $p>1/\liminf_{n \rightarrow \infty}(\theta_{n}{\mathbb E}[T_n])^{1/n}$
while it is 0 if $p<1/\liminf_{n \rightarrow \infty}(\theta_{n}{\mathbb E}[T_n])^{1/n}.$
This gives us that
$p_c(I)=1/\liminf_{n \rightarrow \infty}(\theta_{n}{\mathbb E}[T_n])^{1/n}$ from which it follows
that
${\rm br}I=\liminf_{n \rightarrow \infty}(\theta_{n}{\mathbb E}[T_n])^{1/n}$ (recall that
$p_c(I)=1/{\rm br}I$).

We now proceed with the final step in proving that
$$\liminf_{n \rightarrow \infty}(\theta_{n}{\mathbb E}[T_n])^{1/n}=\liminf_{n \rightarrow
\infty}{\mathbb
E}[T_n]^{1/n}.$$
Obviously, $\theta_{n}{\mathbb E}[T_n] \leq {\mathbb E}[T_n]$ for every $n,$ so
we only need to show that
$$\liminf_{n \rightarrow \infty}(\theta_{n}{\mathbb E}[T_n])^{1/n} \geq \liminf_{n \rightarrow
\infty}{\mathbb
E}[T_n]^{1/n}.$$
As before, let $m=\liminf_{n \rightarrow \infty}{\mathbb E}[T_n]^{1/n}>1$ and
choose $\epsilon>0,$ so that $m (1-\epsilon)>1.$ Furthermore, we can choose an $N$ such that
${\mathbb E}[T_n]^{1/n} \geq m (1-\epsilon)$ for every $n \geq N.$
Using inequalities (\ref{eqn19}) and (\ref{eqn18}) we get that for some constant $C$ and $m \geq n,$
\[
{\mathbb P}(T_m>0 \Big{|} T_n=1) \geq
\left[C{\mathbb E}[T_{n}]\sum_{j=n+1}^{m} \frac{1}{{\mathbb E}[T_{j}]}\right]^{-1}
 \geq \left[C{\mathbb E}[T_{n}]\sum_{j=n+1}^{\infty} \frac{1}{{\mathbb E}[T_{j}]}\right]^{-1}.
\]
Therefore, for $n \geq N,$
\begin{eqnarray*}
\lefteqn{\theta_n {\mathbb E}[T_{n}]= \lim_{m \rightarrow \infty}{\mathbb P}(T_m>0 \Big{|}
T_n=1){\mathbb
E}[T_{n}]} \\
& & \geq \left[C\sum_{j=n+1}^{\infty} \frac{1}{{\mathbb E}[T_{j}]}\right]^{-1}
\geq \left[C\sum_{j=n+1}^{\infty} \frac{1}{(m(1-\epsilon))^{j}}\right]^{-1} \\
& & =\left[\frac{C}{(m(1-\epsilon))^n}\sum_{j=1}^{\infty} \frac{1}{(m(1-\epsilon))^{j}}\right]^{-1}
=(m(1-\epsilon))^n C',
\end{eqnarray*}
where $C'>0.$ Therefore,
for all $n \geq N$ we have
\[
(\theta_n {\mathbb E}[T_{n}])^{1/n} \geq m(1-\epsilon) C'^{1/n},
\]
so that
\[
\liminf_{n \rightarrow \infty}(\theta_n {\mathbb E}[T_{n}])^{1/n} \geq m(1-\epsilon).
\]
Since $\epsilon>0$ can be choosen arbitrarily small, we are done.

\hfill{$\Box$}

\medskip\noindent
{\bf Remark} In fact, the proof of Theorem \ref{thm4} shows that if the family tree $T$ of
a Galton-Watson process satisfies (\ref{eqn11}), (\ref{eqn12}) and (\ref{eqn13}), and
$p>p_c(\overline{T})$ ,
then so
does the family tree associated with the $I^p$ process.

\section{Proof of Theorem \ref{propunifconv} and Proposition \ref{thm2}} \label{sec3}

Before we can prove Theorem \ref{propunifconv}, we need the following domination lemmas.
The first one appears (without proof) in \cite{BHS}. The
proof we give is due to Olle H\"aggstr\"om (unpublished).

\begin{lemma}  \label{lem:conditioned_binomial}
For $k \geq 1$, $p \in (0,1)$ and $0 \leq m \leq k$, write
$\rho_{k,p,m}$ for the distribution of a Binomial$(k,p)$ random variable
conditioned on taking value at least $m$. For $p_1 \leq p_2$, we have
\[
\rho_{k,p_1,m} \, \preceq \, \rho_{k,p_2,m} \, ,
\]
where $\preceq$ denotes stochastic domination.
\end{lemma}

\medskip\noindent
{\bf Proof.} For $i=1,2$, let $Y_i$ be a Bin$(k, p_i)$ random variable, and let
$X_i$ be a random variable with distribution $\rho_{k, p_i, m}$.
Since $x/(1-x) < y/(1-y)$ for $0 < x < y <1$, it is enough to show that for any
$n\in\{m+1, \ldots, k\}$ we have
\[
\frac{{\mathbb P}(X_1 \geq n)}{{\mathbb P}(X_1 < n)} \, \leq \,
\frac{{\mathbb P}(X_2 \geq n)}{{\mathbb P}(X_2 < n)},
\]
which is the same as showing that
\begin{equation}  \label{eq:need_to_show}
\frac{{\mathbb P}(X_2 \geq n)}{{\mathbb P}(X_1 \geq n)} \cdot \frac{{\mathbb P}(X_1 < n)}{{\mathbb
P}(X_2 < n)}
\, \geq \, 1 \, .
\end{equation}
Writing $Z_1$ and $Z_2$ for the probabilities that $Y_1 \geq m$ and
$Y_2 \geq m$, respectively, the left-hand side of (\ref{eq:need_to_show})
becomes
\begin{equation}  \label{eq:first_rewrite}
\frac{\frac{1}{Z_2} \sum_{j=n}^k {k \choose j}
p_2^j(1-p_2)^{k-j}}{\frac{1}{Z_1}
\sum_{j=n}^k {k \choose j} p_1^j(1-p_1)^{k-j}} \cdot
\frac{\frac{1}{Z_1}
\sum_{j=m}^{n-1} {k \choose j} p_1^j(1-p_1)^{k-j}}{\frac{1}{Z_2}
\sum_{j=m}^{n-1} {k \choose j}
p_2^j(1-p_2)^{k-j}} \, .
\end{equation}
Cancelling the $Z_i$'s and introducing the notation
$\phi_i=\frac{p_i}{1-p_i}$ for $i=1,2$, the expression in
(\ref{eq:first_rewrite}) may further be rewritten as
\begin{eqnarray}  \nonumber
\lefteqn{ \mbox{ } \hspace{-20mm}
\frac{p_2^n(1-p_2)^{k-n} \sum_{j=n}^k {k \choose j}
\phi_2^{j-n}}{p_1^n(1-p_1)^{k-n} \sum_{j=n}^k {k \choose j} \phi_1^{j-n}}
\cdot
\frac{p_1^n(1-p_1)^{k-n} \sum_{j=m}^{n-1} {k \choose j}
\phi_1^{j-n}}{p_2^n(1-p_2)^{k-n} \sum_{j=m}^{n-1} {k \choose j}
\phi_2^{j-n}} = } \\
& = &
\frac{ \sum_{j=n}^k {k \choose j}
\phi_2^{j-n}}{ \sum_{j=n}^k {k \choose j} \phi_1^{j-n}}
\cdot
\frac{ \sum_{j=m}^{n-1} {k \choose j}
\phi_1^{j-n}}{ \sum_{j=m}^{n-1} {k \choose j}
\phi_2^{j-n}} \, .
\label{eq:second_rewrite}
\end{eqnarray}
Now note that $\phi_1 \leq \phi_2$, so that
\[
\sum_{j=n}^k {k \choose j}
\phi_2^{j-n} \, \geq \, \sum_{j=n}^k {k \choose j} \phi_1^{j-n}
\]
and
\[
\sum_{j=m}^{n-1} {k \choose j}
\phi_1^{j-n} \, \geq \, \sum_{j=m}^{n-1} {k \choose j}
\phi_2^{j-n} \, .
\]
Hence, the expression in (\ref{eq:second_rewrite}) is greater than or equal
to $1$, so (\ref{eq:need_to_show}) is verified and the lemma is
established. \hfill{$\Box$}
\vspace{3mm}

\noindent
We proceed with the following lemma. We will in fact only use it in the case
$m=1$, but we nevertheless provide a proof of the general statement.

\begin{lemma}
\label{more_bin_cond2}
In the notation of Lemma \ref{lem:conditioned_binomial}, it is the case that for any
$1\leq k \leq l$ and $0 \leq m \leq k$
$$
\rho_{k,p,m} \preceq \rho_{l,p,m},
$$
for all $0 < p < 1.$
\end{lemma}

\medskip\noindent
{\bf Proof.}
It is obvious that we only need to prove the lemma in the case $l=k+1.$ Therefore,
let $Y_1,\ldots, Y_{k+1}$ and $X_1,\ldots, X_k$ be i.i.d.\ Bernoulli random variables with
expectation $p$ and let $Y=\sum_{i=1}^{k+1} Y_i$ and
$X=\sum_{j=1}^k X_j$. We need to show that $\P(X \geq n|X \geq m) \leq \P(Y \geq n|Y \geq m)$, for
all $n=m,m+1,\ldots, k$. To this end we write
\begin{eqnarray*}
\P(Y \geq n|Y\geq m) &=& \P(Y\geq n|Y\geq m, Y_{k+1}=0) \P(Y_{k+1}=0|Y\geq m)\\
& & + \, \P(Y \geq n |Y\geq m, Y_{k+1}=1) \P(Y_{k+1}=1|Y\geq m)\\
&=& \P(X \geq n | X\geq m)\P(Y_{k+1}=0|Y\geq m) \\
& & + \, \P(X \geq n-1| X \geq m-1 )\P(Y_{k+1}=1|Y\geq m).
\end{eqnarray*}
Therefore, we need to show that for $n >m$,
\[
\P(X \geq n-1| X \geq m-1 ) \geq \P(X \geq n| X \geq m),
\]
or equivalently,
\[
\P(X \geq n| X \geq n-1 ) \leq \P(X \geq m| X \geq m-1 ).
\]
It is easy to see that it suffices to prove this for $m=n-1,$ or to simplify notation, to show that
\[
\P(X \geq n+1| X \geq n ) \leq \P(X \geq n| X \geq n-1 ).
\]
Since
\[
\P(X \geq n+1| X \geq n )=1-\P(X =n| X \geq n ),
\]
we need to show that
\[
\P(X =n-1| X \geq n-1 ) \leq \P(X =n| X \geq n ).
\]
Writing $p_n:=\P(X=n)$ we rewrite this as
\[
\frac{p_n+\cdots+p_k}{p_{n-1}+\cdots+p_k} \leq \frac{p_n}{p_{n-1}},
\]
or equivalently that
\begin{equation} \label{eqnnew2}
p_{n-1}(p_{n+1}+\cdots+p_k) \leq p_n (p_{n}+ \cdots+ p_k).
\end{equation}
It suffices to show that $p_{n-1}p_{n+j} \leq p_n p_{n+j-1}$, for $1 \leq j \leq k-n$.
This however is easily checked by a straightforward calculation.
\hfill{$\Box$}
\vspace{3mm}

\noindent
We are now ready to prove Theorem \ref{propunifconv}.

\medskip\noindent
{\bf Proof of Theorem \ref{propunifconv}.} For the purpose of this proof, we introduce a new
stochastic process $\tilde{I}_n^p$, indexed by $n=1,2,\ldots$ as follows. $\tilde{I}_1^p$ is
distributed
as the number of points in $I_1^p$. If this number of points is 0 however, we resample according to
the same distribution, and repeat this until the total number of offspring is at least 1. If we do
{\em
not}
resample at this first generation, we define $R_0:=1$; if we do resample, we set $R_0=0.$

In an inductive fashion, having defined $\tilde{I}_n^p$, we consider all points in $\tilde{I}_n^p$
and
give each of them a random number
of offspring distributed as $L_{n+1}''$, independently of each other. However, if the total number
of
offspring
is 0, we resample {\em all} offsprings using the same distributions, until the total number of
offspring
is at least 1. If we do not have to resample, we define $R_n:=1$; if we do resample, we set $R_n=0$.
Of course, the distribution of the number of points in $\tilde{I}_n^p$ given $\tilde{I}_{n-1}^p=k$
for
some
$k \geq 1$
is the same as the distribution of the number of points in $I_n^p$ given $I_{n-1}^p=k$ conditioned
on
being at
least one.

We can now write, for any $M,$
\begin{eqnarray}
\label{belangrijk}
\lefteqn{\P(0 < I_n^p <M) = \P(\prod_{i=0}^{n-1}R_i =1, 0 < \tilde{I}_n^p < M)} \nonumber\\
& &  \leq \P(0 < \tilde{I}_n^p < M )
= \P(\tilde{I}_n^p < M).
\end{eqnarray}

Now let $p_c(\overline{T}) < p < q$. We claim that
$$
\tilde{I}_n^p \preceq \tilde{I}_n^q.
$$
To see this, we note that the offspring distributions of $I^p$ can be realised by first
drawing from the appropriate $L_n'$, and then keep all points in the offspring with probability
$p$, independently of each other.  Now the combination of Lemma \ref{lem:conditioned_binomial} and
Lemma \ref{more_bin_cond2} implies that for $k \leq \ell$ and $p \leq q$ we have
\begin{equation}
\label{ttrr}
\rho_{k, p, 1} \preceq \rho_{\ell, q, 1}.
\end{equation}
Clearly, we can couple $\tilde{I}_1^p$ and $\tilde{I}_1^q$ so that
$\tilde{I}_1^p \leq \tilde{I}_1^q$, since we can use the same offspring $L_1'$ for
them to get $I_1$
and then the domination follows from Lemma \ref{lem:conditioned_binomial}.
 Let $\{L'_{2,i}\}_{i=1}^{\tilde{I}_1^q}$
be i.i.d.\ with distribution equal to $L'_2$ and independent of everything else.
We can now get $\tilde{I}_2^p$ by letting it be a ${\rm Bin}(\sum_{i=1}^{\tilde{I}_1^p}L'_{2,i},p)$
conditioned on being at least one. Similarly we get
$\tilde{I}_2^q$ by letting it be a ${\rm Bin}(\sum_{i=1}^{\tilde{I}_1^q}L'_{2,i},q)$
conditioned on being at least one.
The fact that we can couple $\tilde{I}_2^q$ and $\tilde{I}_2^p$ so that
$\tilde{I}_2^q \leq \tilde{I}_2^p$ now
follows from (\ref{ttrr}).
Repeating this procedure at every level gives that
\begin{equation} \label{eqn20}
\P(\tilde{I}_n^p<M) \leq \P(\tilde{I}_n^{p_1}<M),
\end{equation}
for all $p > p_1$, and this is where the uniformity in $p$ comes from.

Of course letting $M$ above depend on $n$ does not change the validity of the argument.
According to (\ref{belangrijk}) and (\ref{eqn20}) it therefore
suffices to show that
$$
\P(\tilde{I}_n^{p_1} < ((1-\epsilon){\rm br} \overline{I^{p_1}})^n) \to 0,
$$
as $ n \to \infty.$ For this, we use Theorem \ref{thm4} and Proposition
\ref{lemma5}. Consider the subsequence $\{n_k\}$ and the constant $C>0$ dictated by applying
Proposition
\ref{lemma5} to $I^{p_1}$. This
is allowed according to the remark following the proof of Theorem \ref{thm4}. Since each element in
the
$n_1$th generation of the $I^{p_1}$ process has a probability at least $C$ to survive, there
is at least probability $C>0$ that no resampling is ever going to be necessary in the
$\tilde{I}^{p_1}$
process after time $n_1$. There are now two possibilities. Either, at some point resampling
is needed, or no resampling is ever needed after time $n_1$.

In the latter case, we have that $\tilde{I}_n^{p_1}$ is at least as large as the number of
points in a surviving copy of an $I^{p_1}$ tree with only one vertex at generation $n_1.$ It follows
from
Theorem \ref{thm4} that this surviving tree has branching number ${\rm br} \overline{I^{p_1}}.$
Using that the lower growth number is at least as large as the branching number we are done
in this case.

On the other hand, if resampling is needed, then we take the first element in the subsequence
$\{n_k\}$ after the first resampling, and repeat the reasoning from there. It follows that a.s.,
$\liminf_{n \rightarrow \infty} (\tilde{I}_n^{p_1})^{1/n} \geq {\rm br} \overline{I^{p_1}}$,
and the proof is complete. \hfill{$\Box$}
\vspace{3mm}

\noindent
We can now prove Proposition \ref{thm2}

\medskip\noindent
{\bf Proof of Proposition \ref{thm2}}. We write
$$
\theta(p)=\P(I_n^p>0)-\P(I_n^p>0, |I^p| < \infty),
$$
recall that $I^p$ denotes the component of the root. We will prove that along a subsequence,
the last term tends to zero uniformly in $p \in [p_1, 1]$, where $p_1 > p_c(\overline{T})$,
from which the result follows.

Since the $p$-dependence is important now, we write $\theta_n(p)$ for $\theta_n$ in the context
of the Galton-Watson process associated with $I^p$. For any $M >0$ we  write, for $p_1 \leq p \leq
1$,
\begin{eqnarray*}
\P(I_n^p>0, |I^p| < \infty)
& \leq  & \P(0 < I_n^p < M) + \P(I_n^p \geq M, |I^p| < \infty) \\
& \leq & \P(0 < I_n^p < M) + (1-\theta_n(p))^M \\
& \leq &  \P(0 < I_n^p < M) + (1-\theta_n(p_1))^M.
\end{eqnarray*}

\noindent
Let $\epsilon >0$ be arbitrary. We want to apply Proposition \ref{lemma5} to $I^p.$ According to the
remark after the proof of Theorem \ref{thm4}, all the assumptions of Proposition \ref{lemma5}
holds for $I^p$ since $p>p_c(\overline{T}).$

Now let $C$ be the constant in Proposition \ref{lemma5} when we apply it to $I^{p_1}$. We choose $M$
so
large that $(1-C)^M < \epsilon/2$. Next choose $n$ in the appropriate
subsequence of Proposition \ref{lemma5} and at the same time so large that
the first term at the right hand side is at most $\epsilon/2$; this is possible according to
Theorem \ref{propunifconv} above. The right hand side is then bounded
above by $\epsilon$, uniformly in $p \in [p_1,1]$.
In summary, for any $\epsilon >0$ we can find $K$ such that
$$
\theta(p) \geq \P(I_{n_k}^p>0 )-\epsilon
$$
for every $p \in [p_1,1]$ and every $n_k$ in the subsequence dictated by Proposition \ref{lemma5}
with
$k
\geq K.$
We see that for all $k \geq K$ and for all $p \in [p_1,1]$,
\[
|\theta(p)-\P(I_{n_k}^p>0)| \leq \epsilon,
\]
which concludes the argument.
\hfill{$\Box$}
\vspace{3mm}

\section{Continuity of the percolation function}
\label{sec4}

The supercritical continuity of $\theta(p)$ (Corollary \ref{corr1})
follows immediately from Proposition \ref{thm2}. We point out however that it is possible to obtain
the same result by combining Theorem \ref{thm4} with a modified version of the classical
argument found in \cite{BK}. We provide a sketch.

\medskip\noindent
{\bf Sketch of proof of Corollary \ref{corr1} from Theorem \ref{thm4}.}
We start by drawing an $I$ from the correct distribution. Associate to every
edge $e$ in $I$ an independent $U([0,1])$ random variable, denoted by $U_e.$
For $p_c<q<p,$ create ${\cal I}^q$ and ${\cal I}^p$ by keeping every vertex of $I$ and
those edges $e \in I$ such that $U_e \leq q,p$ respectively. Consider any infinite subtree
$J$ in ${\cal I}^p.$ Theorem \ref{thm4} gives us that $p_c(J)=1/{\rm br}J=1/p\liminf_{n \to \infty}
{\mathbb E}[I_n]^{1/n}=p_c(I)/p$ a.s. Therefore, performing further percolation on $J$ with
density $q/p>p_c(I)/p$ will result in a new graph containing an infinite
subgraph a.s. Of course, the distribution of this new graph must be the same as
$J \cap {\cal I}^q.$ Furthermore this holds in particular if $J=I^p$ showing that
if $|I^p|=\infty$ then there exists a.s. an infinite subtree of $I^p \cap {\cal I}^q.$ It is now
possible to proceed as in \cite{BK}. \hfill{$\Box$}

\medskip
The non-classical way to conclude continuity of the percolation function has an interesting analogy
on
$\mathbb{Z}^d$. Define $B_n:= [-n,n]^d$ and write $\partial B_n$ for the (inner) boundary of $B_n.$
Letting $\{0 \leftrightarrow \partial B_n\}$ denote the event that the origin is connected
to $\partial B_n$ by a path of open edges,
define
\[
\varphi_n(p):={\mathbb P}_p(0 \leftrightarrow \partial B_n).
\]
Clearly,
\begin{equation}
\label{convergence}
\theta(p)=\lim_{n \to \infty} \varphi_n(p),
\end{equation}
for all $0 \leq p \leq 1$. The inequality of the following equation (valid for every $n \geq 1$)
is a part of Theorem 8.18 of \cite{G}:
\begin{equation} \label{eqn22}
\varphi_n(p)-\theta(p)=\P_p(0 \leftrightarrow B_n, |C| <\infty) \leq A(p,d)n^de^{-n \sigma(p)},
\end{equation}
where we can take
\begin{equation} \label{eqn23}
A(p,d)=\frac{d^2}{p^2(1-p)^{d-2}}.
\end{equation}
Furthermore, according to Theorem 8.21 of \cite{G} we can take $\sigma(p)$ to be uniformly bounded
away from $0$ on any closed sub-interval of $(p_c,1).$ We point out the following
corollary and sketch the proof.
\begin{corollary}
The percolation function $\theta(p)$ on ${\mathbb Z}^d,$ $d \geq 2$ is continuous
for $p>p_c.$
\end{corollary}
{\bf Sketch of proof.} Choose $p_c<p_1<p_2<1.$
Combining equations (\ref{eqn22}), (\ref{eqn23}) and Theorem 8.21 of \cite{G}
explained directly above, it is straightforward to prove that there exists constants
$C=C(p_1,p_2)<\infty$ and $\delta=\delta(p_1,p_2)>0$ such that for any $p \in [p_1,p_2]$
and any $n \geq 1,$
\[
\varphi_n(p)-\theta(p) \leq Ce^{-n \delta}.
\]
Since trivially
\[
\theta(p) \leq \varphi_n(p),
\]
it follows that $\varphi_n(p) \to \theta(p)$ uniformly on any closed subinterval of $(p_c,1),$
from which the statement follows.
\hfill{$\Box$}


\vspace{3mm}









\begin{thebibliography}{99}

\bibitem{Agresti} Agresti, A.\ {\em On the extinction times of varying and random environment
branching processe}, J.\ Appl.\ Prob.\ {\bf 12}, 39-46 (1975).

\bibitem{BK} Van den Berg, J.\ and Keane, M., {\em On the continuity of the percolation probability
function}, Conference in modern analysis and probability, R.\ Beats et al
(ed),
61 - 65  AMS, Providence, RI (1982).

\bibitem{BHS} Broman E.I., H\"{a}ggstr\"{o}m O. and Steif J. E.,
{\em Refinements of Stochastic Domination},
Probab. Theory and Rel.\ Fields {\bf 136} No. 4, 587-603 (2006).

\bibitem{G} Grimmett, G., {\em Percolation},
Second edition, Springer-Verlag, Berlin (1999).

\bibitem{HP} H\"aggstr\"om, O.\ and Peres, Y., {\em Monotonicity of uniqueness for percolation on
Cayley graphs: all infinite clusters are born simultaneously}, Probab. Theory and Rel. Fields {\bf
113},
273 - 285 (1999).

\bibitem{Lyons90} Lyons, R., {\em Random walks and percolation on trees},
 Ann. Probab. {\bf 18}  no. 3, 931--958 (1990).

\bibitem{Ly2} Lyons, R., {\em Random walks, capacity and
percolation on trees},  Ann.\ Prob.\ {\bf 20} no. 4, 2043-2088 (1992).

\bibitem{R} Lyons, R., {\em Probability on trees and networks}, In progress,
URL: http://mypage.iu.edu/~rdlyons/prbtree/prbtree.html.

\bibitem{S} Schonmann, R., {\em Stability of infinite clusters in supercritical percolation},
Probab. Theory and Rel. Fields  {\bf 113}, 287 - 300 (1999).

\end{thebibliography}
\end{document}